\documentclass[12pt]{amsart}
\usepackage{psfig}
\usepackage{multicol}
\newcommand{\C}{\mathbb C}
\newcommand{\R}{\mathbb R}
\newcommand{\Z}{\mathbb Z}

\newcommand{\dd}{\partial}

\newtheorem{thm}{Theorem}

\newtheorem{prop}{Proposition}
\newtheorem{cor}{Corollary}

\theoremstyle{definition}

\theoremstyle{remark}

\newtheorem{rmk}{Remark}

\newcommand{\conj}{\operatorname{conj}}
\newcommand{\indx}{\operatorname{ind}}
\begin{document}

\title[Adjunction inequality for real curves]
{Adjunction inequality for real algebraic curves}
\author{G.Mikhalkin}

\thanks{Research at MSRI supported in part by NSF grant \#DMS 9022140}

\address{Dept. of Mathematics\\
University of Toronto\\ Toronto, Ont. M5S 3G3\\ Canada}
\curraddr{MSRI\\ 1000 Centennial Dr.\\ Berkeley, CA 94720\\ USA}

\email{mihalkin@math.toronto.edu, grisha@msri.org}
\begin{abstract}
The zero set of a real polynomial in two variable
is a curve in $\R^2$.
For a generic choice of its coefficients this is a non-singular
curve, a collection of circles and lines properly embedded in $\R^2$.
What topological arrangements of these circles and lines
appear for the polynomials of a given degree?
This question arised in the 19th century in the works
of Harnack and Hilbert and was included by Hilbert into
his 16th problem.
Several partial results were obtained since then
(see \cite{Pe}, \cite{A}, \cite{R}, \cite{V}).
However the complete answer is known only
for polynomials of degree 5 (see \cite{Po}) or less.
The paper presents a new partial result toward
the solution of the 16th Hilbert problem.

A consequence of the main theorem of the paper
implies that the zero set of a polynomial of degree 7
can not be isotopic to a curve pictured on Fig.\ref{deg7}
(contrary to the curve pictured on Fig.\ref{deg7y}).
The realizibility of Fig.\ref{deg7} does not
contradict to the previously known restrictions.
The proof of the main theorem makes use of the proof
by Kronheimer and Mrowka \cite{KM} of
the Thom conjecture in $\C P^2$.
\end{abstract}
\maketitle

\begin{figure}[h]
\begin{multicols}{2}
\centerline{
\psfig{figure=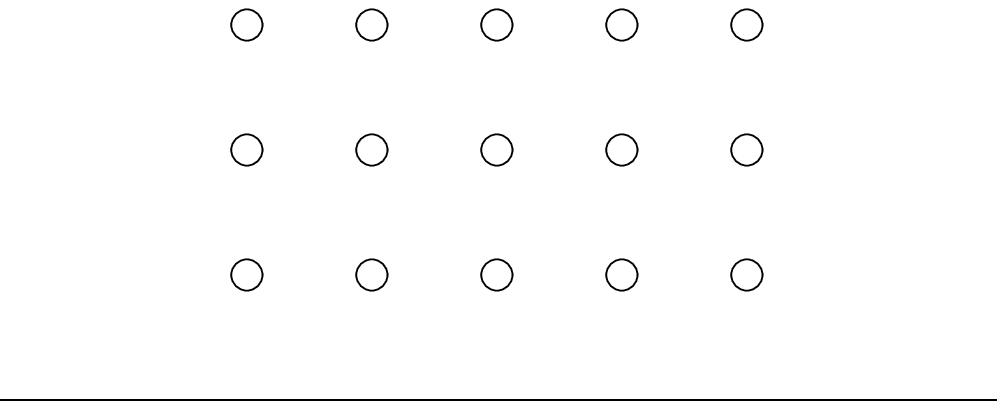,height=1in,width=2in}}
\caption{\label{deg7} This picture is not isotopic
to a curve of degree 7.}
%\end{figure}
%\begin{figure}[h]
\centerline{
\psfig{figure=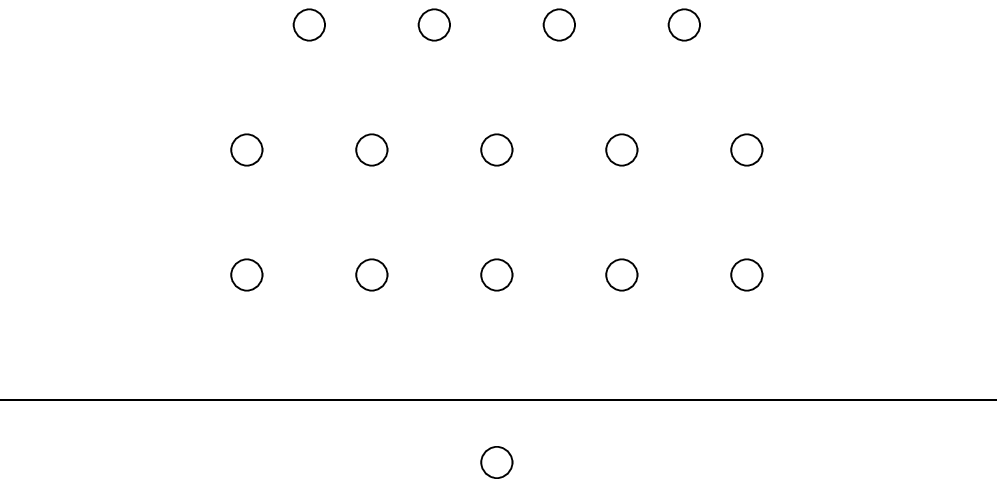,height=1in,width=2in}}
\caption{\label{deg7y} This picture is isotopic
to a curve of degree 7.}
\end{multicols}
\end{figure}
\section{Statement of the theorem}
Let $f$ and $g$ be two real homogeneous polynomials in
two variables of degrees $m$ and $n$, $m\ge n$, $m\equiv n\pmod{2}$.
Let
\begin{gather*}
\R A=\{(x_0:x_1:x_2)\in\R P^2\ |\ f(x_0,x_1,x_2)=0\}\\
\R B=\{(x_0:x_1:x_2)\in\R P^2\ |\ g(x_0,x_1,x_2)=0\}
\end{gather*}
be the corresponding curves in $\R P^2$.
Let
\begin{gather*}
\C A=\{(x_0:x_1:x_2)\in\C P^2\ |\ f(x_0,x_1,x_2)=0\}\\
\C B=\{(x_0:x_1:x_2)\in\C P^2\ |\ g(x_0,x_1,x_2)=0\}
\end{gather*}
be their complexifications in $\C P^2$.
If the coefficients of $f$ and $g$ are generic then
$\C A$ and $\C B$ are transverse smooth submanifolds of $\C P^2$.
We call such curves {\em non-singular}.
The real parts $\R A$ and $\R B$ of non-singular curves
are transverse smooth submanifolds of $\R P^2$.

The complex conjugation involution
$\conj:(x_0:x_1:x_2)\to (\bar{x}_0:\bar{x}_1:\bar{x}_2)$
acts on $\C P^2$, $\C A$ and $\C B$
fixing $\R P^2$, $\R A$ and $\R B$.
Note that $\C A-\R A$ is either conneceted (then $\R A$
is called {\em of type II}) or consists of two components
interchanged by $\conj$ (then $\R A$ is called {\em of
type I}) since $\C A/\conj$ is connected and $\dd\C A/\conj=\R A$.
If $\R A$ is of type I then following \cite{R}
we equip $\R A$ with {\em complex orientation},
i.e. with the boundary orientation induced from
one of the components of $\C A-\R A$.
Because of ambiguity in the choice of the component of $\C A-\R A$
this orientation is defined up to simultaneous reversion of
the orientations of all the components of $\R A$.

Let $C\subset\R P^2$ be an oriented null-homologous
topologically embedded (multicomponent) curve.
Following \cite{V} we define a locally constant
function called the {\em index}
$$|\indx|:\ \R P^2-C\ \to\ \Z$$ in the following way.
The index of the non-orientable component of $\R P^2-C$ is zero.
The complement of the non-orientable component is
a disjoint union of disks.
A regular neighborhood $N$ of each of these disks is $\R^2$
and in $N$ there is a well-defined index of a point with
respect to $C\cap N$ once we pick an orientation of $N$.
We define $\indx|_{N-C}$ as the absolute value of the index in $N$.
\begin{figure}[h]
\centerline{\psfig{figure=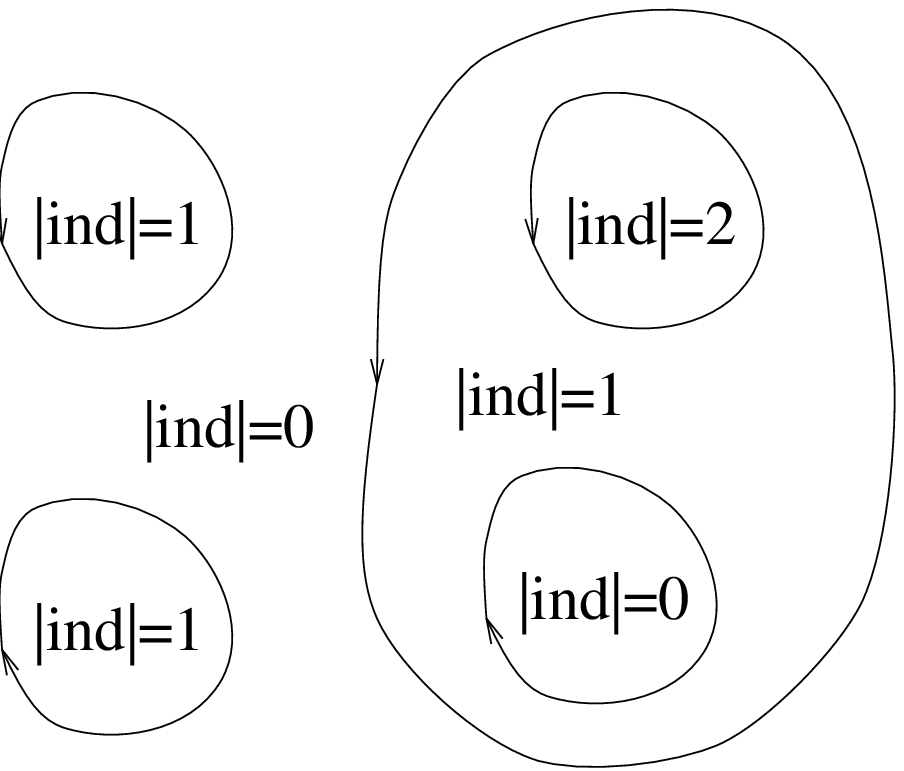,height=1in,width=1.5in}}
\caption{\label{index} Function $|\indx|$.}
\end{figure}

For each locally constant function $F:\R P^2-C\to\R$ we define
$$\int\limits_{\R P^2}Fd\chi=\sum_G F(G)\chi(G)$$
over all components $G$ of $\R P^2-C$
($\chi$ denotes the Euler characteristic).

\begin{figure}[h]
\centerline{\psfig{figure=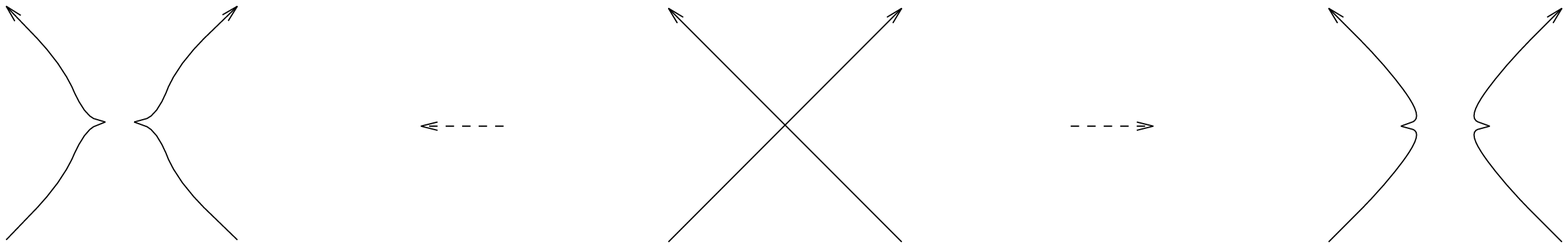,height=0.66in,width=4.7in}}
\caption{\label{sm12} Smoothings of a double point.}
\end{figure}

Let $\R A$ and $\R B$ be the curves of type I equipped with
the complex orientations.
We associate the following curve $C$ to $\R A$ and $\R B$.
Take $\R A\cup\R B$ and smooth all the points of
intersection according to their orientations
choosing any pattern of Fig.\ref{sm12} (cf. smoothing in \cite{MP}).
(Recall that $\R A$ and $\R B$ are transverse
since we assume that coefficients of $f$ and $g$ are generic
and $C$ is homologous to zero since $m\equiv n\pmod{2}$).

Let $G$ be a domain in $\R P^2$ such that $\dd G$
is a component of $C$.
Let $c_{\operatorname{out}}$ ($c_{\operatorname{in}}$)
be the number of the outward (inward) with respect to $G$
cusps on the boundary of $G$.
Let $\chi'(G)=\chi(G)-\frac{1}{2}(c_{\operatorname{out}}-
c_{\operatorname{in}})$ (cf. Remark 1.3 of \cite{MP}).
Define "the negative variation" of $\chi'$ by the following rules.
For a connected $G$
let $\chi'_-(G)=-\chi'(G)$ if $\chi'(G)<0$ and
$\chi'_-(G)=0$ if $\chi'(G)\ge 0$.
For a disconected $G$ define $\chi'_-(G)$ as the sum of $\chi'$
over all components of $G$.
Define
$$\Xi=
\sum_{s=1}^{\infty}s\chi'_-\{|\indx(x)|>s\}.$$
Note that $\Xi=0$ if $|\indx|<2$ on all components of $\R P^2-C$
of negative Euler characteristic.
This will be the case for our applications.
Let $d$ be the number of the intersection points of $\R A$ and $\R B$.

\begin{thm}
$$\int\limits_{\R P^2}|\indx|d\chi\le\frac{(m+n)^2}{4}-3n+d+\Xi.$$
\end{thm}
Note that there is ambiguity in choosing the complex orientations
of $\R A$ and $\R B$ and in smoothing of double points of $\R A\cup\R B$.
The inequality of the theorem holds for {\em any} choice of
the complex orientation and smoothings.

\begin{rmk}
If $\R A$ is a non-singular curve of type I of degree $m$
then the curve obtained by taking $s$ parallel copies of
each component of $\R A$ is isotopic to a non-singular
curve of type I of degree $sm$.
Indeed, if $f$ defines $\R A$ then
$$f(f-\epsilon h)\dots(f-(s-1)\epsilon h)=\delta p$$
defines the requred curve for a choice of real numbers
$\epsilon$, $\delta$ and polynomials $h$ and $p$ of degree $m$ and $sm$.

This observation combined with the theorem yields additional
restrictions on mutual position of two curves of type I.
In particular, it allows to use the theorem in the case when
the degrees of two curves are of opposite parity.
\end{rmk}

\section{Proof of Theorem 1}

Consider first the special case when $\R A\cap\R B=\emptyset$.
Let $A_+$ be the component of $\C A-\R A$ which determines
the chosen complex orientation of $\R A$ and let $B_-$
be the component of $\C B-\R B$ which determines the orientation
{\em opposite} to the chosen orientation of $\R B$.
Let
$$F'=A_+\cup B_-\cup L,$$
where $L=\bigcup_{s=1}^{\infty}\{|\indx|\ge s\}$.
Note that $F'$ makes an integer cycle such that the orientation
induced from $F'$ on $A_+$ is holomorphic and the orientation
induced from $F'$ on $B_-$ is {\em anti}holomorphic.

Let $F$ be the surface smoothly generically immersed into $\C P^2$
obtained by smoothing the corners and, then, perturbation of $F'$.
The double points of $F$ come from two sources, from $A_+\cap B_-$
and from perturbation of $\bigcup_{s=1}^{\infty}\{|\indx|\ge s\}$.
Note that the double points coming from $A_+$ are all negative
since all points of intersection of $\C A$ and $\C B$ equipped
with the holomorphic orientation are positive.

To obtain $F$ out of $F'$ we perturb every component $G'\subset\R P^2$
of $\{|\indx|\ge s\}$, $s>1$, to $G$ along a generic vector field $v$
normal to $G'$ but tangent to $\C A\cup\C B$ along
$\dd G'\subset\R A\cup\R B$ (to keep $F$ a cycle).
The field $iv$ is tangent to $G'$ and to $\dd G'$.
The zeroes of $v$ and $iv$ are of opposite index.
Therefore (inductively by $s$)
we can choose the perturbation field $v$
so that $G$ (equipped with the orientation induced from $F$)
intersects the result of perturbation of
$\{|\indx|\ge t\}$, $t<s$, in $|\chi(G')|$ points,
which are negative if $\chi(G')>0$
and positive if $\chi(G')<0$.
Therefore, the total number of positive double points of $F$
is $\Xi$.

Note that
$$[F]=\frac{m-n}{2}[\C P^1]\in H_2(\C P^2).$$
since $F\cup\conj(F)$ is homologous to
$[\C A]-[\C B]=(m-n)[\C P^1]$
and $\conj$ acts as $(-1)$ on $H_2(\C P^2)$.
We need the following strengthening of
the Thom conjecture \cite{KM}
\begin{equation}
\label{tc}
-\chi(F)+2\Xi\ge\frac{(m-n)^2}{4}-3\frac{m-n}{2}.
\end{equation}
This strengthening is implicitly contained in \cite{KM}.
Indeed, it is proven there that if $H$ is a surface
smoothly embedded to $\C P^2\#\bar{\C P}^2\#\dots\#\bar{\C P}^2$
and $H.H\ge 0$ then $-\chi(H)\ge H.H+K.H$, where
$K=(-3,1,\dots,1)\in H_2(\C P^2\#\bar{\C P}^2\#\dots\#\bar{\C P}^2)$
is the canonical class.
We can get $H$ out of $F$ by blowing up all the double points of $F$.
The intersection of $F$ with a small sphere centered at
a double point of $F$ is the Hopf link
(the difference between the positive and negative
Hopf links disappears if we forget the orientation of the links).
The resulting $H$ intersects each exceptional divisor at
two points, which are of the same sign (which is positive)
for positive double points
and of the opposite signs for negative double points.
Therefore, $H.H=F.F-4\Xi$ and $K.H=-3\frac{m-n}{2}+2\Xi$
while $\chi(H)=\chi(F)$.
This implies \eqref{tc}.
Note that the assumption $m\ge n$ insures $H.[\omega]>0$ and
thus enables us to apply \cite{KM}.

By additivity of $\chi$
\begin{multline}
\label{chiF}
\chi(F)=\chi(A_+)+\chi(B_-)+\chi(L)=\\
\chi(A_+)+\chi(B_-)+\sum_{s=1}^{\infty}\chi(|\indx|\ge s)=
\chi(A_+)+\chi(B_-)+\int\limits_{\R P^2}|\indx|d\chi.
\end{multline}
But $\chi(A_+)=\frac{1}{2}\chi(\C A)=-\frac{m^2-3m}{2}$
and $\chi(B_-)=\frac{1}{2}\chi(\C B)=-\frac{n^2-3n}{2}$
by the adjunction formula for the holomorphic curves $\C A$ and $\C B$.
Combining this with \eqref{tc} and \eqref{chiF} we get
$$\int\limits_{\R P^2}|\indx|d\chi\le\frac{m^2-3m}{2}+\frac{n^2-3n}{2}
-\frac{(m-n)^2}{4}+3\frac{m-n}{2}+2\Xi=
\frac{(m+n)^2}{4}-3n+2\Xi.$$

\begin{figure}
\centerline{\psfig{figure=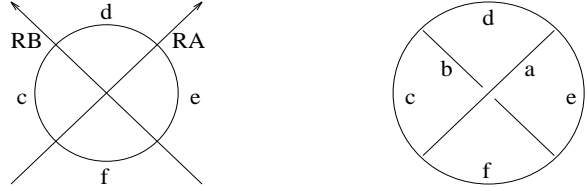,height=1in,width=3in}}
\caption{\label{link} A double point in the real plane
and its link in the complexification.}
\end{figure}

\begin{figure}[b]
\centerline{\psfig{figure=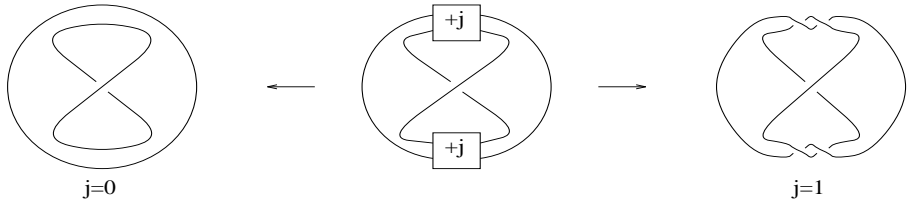,height=1in,width=4.7in}}
\caption{\label{smoo} Choosing the framing along $d$ and $f$.}
\end{figure}

We return to the general case when
$\R A$ and $\R B$ intersect in $d$ points.
We do the same thing still starting from
$F'=A_+\cup B_-\cup\bigcup_{s=1}^{\infty}\{|\indx|\ge s\}$
but we have to pay special attention to smoothing of $F'$
near $\R A\cap\R B$.
The intersection of $A_+$, $B_-$ and $\R P^2$ with
a small sphere centered in a point of $\R A\cap\R B$
is pictured on Fig.\ref{link}.
It consists of 6 arcs,
2 of them, $a$ and $b$ come from the intersection with $A_+$ and $B_-$
and 4 of them, $c$, $d$, $e$, $f$,
come from the intersection with $\R P^2$.
For perturbation of $F'$ we need to choose such a framing
on the circle $a\cup d\cup b\cup f$ that the result of
perturbation of $a\cup d\cup b\cup f$ along this framing
is not linked with $c\cup d\cup e\cup f$.
The number of (full right) twists $j$
of the framing along $f$ has to be equal
to the number of twists along $d$ in order for these circles
to be {\em homologically} unlinked.
These circles are not linked if (and only if) $j=0\ \text{or}\ 1$
(see Fig.\ref{smoo}).
This corresponds to the framing obtained by multiplication by $i$
of the vector field tangent to the result of smoothing with
two inwards or outwards cusps (cf. Fig.\ref{sm12}).
Gluing the disk bounded by $a\cup d\cup b\cup f$
along the arcs $f$ and $d$ decreases the Euler characteristic
of $F$ by 1, so the number of the points of $\R A\cap\R B$
appears in the final formula.

\section{Applications}
By the Harnack inequality \cite{H} an affine curve of degree $2k+1$
has no more then $k(2k-1)$ ovals.
\begin{prop}
If an affine curve of degree $2k+1$ consists of an arc
and $k(2k-1)$ ovals with disjoint interiors then
each component of the complement of the arc in the affine plane
contains at least $\frac{(k-1)(k-2)}{2}$ ovals.
\end{prop}
\begin{rmk}
This proposition is sharp in the following sense.
For any $k>0$ there exist affine curves of degree $(2k+1)$
which consist of $k(2k-1)$ ovals and an arc
where the arc separates $\frac{(k-1)(k-2)}{2}$ ovals from
$\frac{(3k-2)(k+1)}{2}$ ovals.
Construction of such curves can be easily extracted from
Harnack's construction of curves with $(k(2k-1)$ ovals \cite{H}.
\end{rmk}
\begin{proof}
We apply Theorem 1 to the union of the affine curve and
the infinite line in $\R P^2$.
The projectivization of the affine curve is of type I
because the curve is extremal for the Harnack inequality (see \cite{R}).
Suppose the arc separates $A$ ovals from $k(2k-1)-A$ ovals.
Then for a proper choice of the orientation of the infinite
line and the affine curve
$$\int\limits_{\R P^2}|\indx|d\chi\ge k(2k-1)-A + 1-A$$
(the right hand side is the Euler characteristic of
the region with $|\indx|=1$, the only possible regions with
$|\indx|>1$ are the interiors of some of the $A$ ovals
(this depends on their orientation), $d=1$ and $\Xi=0$.
Theorem 1 implies now that
$1+k(2k-1)-2A\le (k+1)^2-3+1$ and, therefore
$A\ge\frac{(k-1)(k-2)}{2}$.
\end{proof}
\begin{cor}
The curve pictured on Fig.\ref{deg7} is not isotopic
to an affine algebraic curve of degree 7.
\end{cor}
\begin{rmk}
Degree 7 is the smallest degree when such a restriction appears.
For an affine curve of degree 5 which consists of an arc and
some ovals there are no constraints on distribution of the ovals
between the half-planes (as long as their total number is no more
than 6).
\end{rmk}

\begin{figure}[h]
\centerline{\psfig{figure=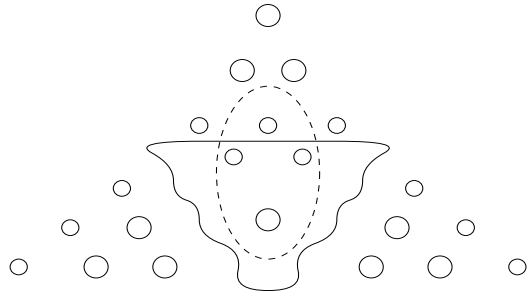,height=1.5in,width=2.7in}}
\caption{\label{deg8} The Harnack curve of degree 8 and a conic.}
\end{figure}

For an even degree $2k$ Harnack's construction
(for the maximal number of ovals) produces
a curve $\R P^2$ which consists of a non-empty oval,
$\frac{k^2-3k+2}{2}$ empty ovals inside the non-empty ovals
and $\frac{3k^2-3k}{2}$ empty ovals outside of the non-empty oval
(see Fig.\ref{deg8} for $k=4$).
Let $\R A\subset\R P^2$ be a curve of degree $2k$ with such
an arrangement of the ovals.
\begin{prop}
\label{ev}
Not more than $(k-2)$ out of the $\frac{k^2-3k+2}{2}$ interior ovals
of $\R A$ may be contained inside an ellipse (a non-singular curve
of degree 2) in $\R P^2$ not intersecting $\R A$.
\end{prop}
\begin{rmk}
This proposition is also sharp.
For the Harnack curves of degree $2k$ we may find ellipses enclosing
$(k-2)$ of the interior ovals and not intersecting the curves.
\end{rmk}
\begin{figure}[h]
\centerline{\psfig{figure=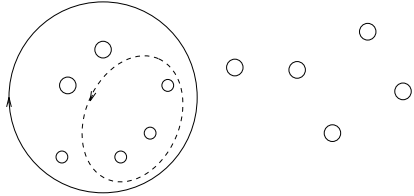,height=1in,width=2.1in}}
\caption{\label{prop2} Orientations of the ellipse and
the non-empty oval of $\R A$.}
\end{figure}
\begin{proof}
Let $L>0$ be the number of the interior ovals inside the ellipse.
Apply Theorem 1 to $\R A$ and the ellipse oriented so that
the non-empty oval of $\R A$ and the ellipse are oriented
in the opposite way (see Fig.\ref{prop2}).
We get
$$\int\limits_{\R P^2}|\indx|d\chi\ge\frac{3k^2-3k}{2}-\frac{k^2-3k+2}{2}+2L=
k^2-1+2L,$$
$d=0$ and $\Xi=0$.
Therefore, $k^2-1+2L\le (k+1)^2-6$ and $L\le k-2$.
\end{proof}
\begin{cor}
It is not possible to separate the 3 interior ovals of a curve of degree 8
with topological arrangement $<18\sqcup 1\!<\! 3\!>>$ from the rest
of the ovals by a conic.
\end{cor}

\begin{rmk}
Theorem 1 generalizes for curves
on algebraic surfaces other than $\R P^2$.
Instead of $m\equiv n\pmod{2}$ we require that
the union of the two curves is null-homologous
in the surface (so that $|\indx|$ is well-defined).
The formula then turns to
$$\int\limits_{\R P^2}|\indx|d\chi\le\frac{([\C A]+[\C B])^2}{4}+
K.[\C B]+d+\Xi,$$
where $K$ is the canonical class of the complexification of the surface.

If the geometric genus of the surface is positive then
we do not need to require that the degree of $\C A$
is greater than or equal to the degree of $\C B$.
\end{rmk}

\end{document}